\documentclass[12pt]{amsart}
\begin{document}

\def\1#1{\overline{#1}}
\def\2#1{\widetilde{#1}}
\def\3#1{\widehat{#1}}
\def\4#1{\mathbb{#1}}
\def\5#1{\frak{#1}}
\def\6#1{{\mathcal{#1}}}
\def\Cal#1{{\mathcal{#1}}}

\def\C{{\4C}}
\def\R{{\4R}}
\def\N{{\4N}}
\def\Z{{\4Z}}

\title[Extension of CR-functions]{Extension of CR-functions on wedges}
\author[D. Zaitsev and G. Zampieri]{Dmitri Zaitsev and Giuseppe Zampieri}
\address{D. Zaitsev: Mathematisches Institut, Eberhard-Karls-Universit\"at T\"ubingen, 72076 T\"ubingen, Germany}
\email{dmitri.zaitsev@uni-tuebingen.de}
\address{G. Zampieri: Dipartimento di Matematica, Universit\`a di Padova, via Belzoni 7, 35131 Padova, Italy}
\email{zampieri@math.unipd.it}
%\subjclass{}

\maketitle
%\tableofcontents

%\def\Label#1{\label{#1}{\bf (#1)}~}
\def\Label#1{\label{#1}}

% Standard sets

\def\cn{{\C^n}}
\def\cnn{{\C^{n'}}}
\def\ocn{\2{\C^n}}
\def\ocnn{\2{\C^{n'}}}

% Abbreviations
\def\dist{{\rm dist}}
\def\const{{\rm const}}
\def\rk{{\rm rank\,}}
\def\id{{\sf id}}
\def\aut{{\sf aut}}
\def\Aut{{\sf Aut}}
\def\CR{{\rm CR}}
\def\GL{{\sf GL}}
\def\Re{{\sf Re}\,}
\def\Im{{\sf Im}\,}

\def\codim{{\rm codim}}
\def\crd{\dim_{{\rm CR}}}
\def\crc{{\rm codim_{CR}}}

\def\phi{\varphi}
\def\eps{\varepsilon}
\def\d{\partial}
\def\a{\alpha}
\def\b{\beta}
\def\g{\gamma}
\def\G{\Gamma}
\def\D{\Delta}
\def\Om{\Omega}
\def\k{\kappa}
\def\l{\lambda}
\def\L{\Lambda}
\def\z{{\bar z}}
\def\w{{\bar w}}
\def\Z{{\1Z}}
\def\t{\tau}
\def\th{\theta}

\emergencystretch15pt
\frenchspacing

\newtheorem{Thm}{Theorem}[section]
\newtheorem{Cor}[Thm]{Corollary}
\newtheorem{Pro}[Thm]{Proposition}
\newtheorem{Lem}[Thm]{Lemma}

\theoremstyle{definition}\newtheorem{Def}[Thm]{Definition}

\theoremstyle{remark}
\newtheorem{Rem}[Thm]{Remark}
\newtheorem{Exa}[Thm]{Example}
\newtheorem{Exs}[Thm]{Examples}

\def\bl{\begin{Lem}}
\def\el{\end{Lem}}
\def\bp{\begin{Pro}}
\def\ep{\end{Pro}}
\def\bt{\begin{Thm}}
\def\et{\end{Thm}}
\def\bc{\begin{Cor}}
\def\ec{\end{Cor}}
\def\bd{\begin{Def}}
\def\ed{\end{Def}}
\def\br{\begin{Rem}}
\def\er{\end{Rem}}
\def\be{\begin{Exa}}
\def\ee{\end{Exa}}
\def\bpf{\begin{proof}}
\def\epf{\end{proof}}
\def\ben{\begin{enumerate}}
\def\een{\end{enumerate}}

\section{Introduction}

The celebrated {\sc Boggess-Polking} theorem \cite{BP82} 
extending classical results of {\sc Lewy} \cite{L56}
states that all (continuous) CR-functions on a generic submanifold $M$ in $\C^N$
extend holomorphically to a wedge in the direction of the convex cone
spanned by the values of the (vector-valued) Levi form of $M$.
Here one starts with a submanifold $M$ and ends with a wedge.
A natural question is to obtain a result generalizing the above one within the category of wedges
with description of the additional directions of extendibility
(see {\sc Tumanov} \cite{T95a,T95b} for general extension results 
to wedges without prescription of the direction).

Given a wedge $V$
in a submanifold $M$ of $\C^N$ (see \S\ref{basic} for the definition) 
whose edge $E\subset \d V$ at a point $p\in E$ is a generic
(i.e.\ $T_pE+iT_pE=T_p\C^N$), in general,
CR-functions on $V$ do not extend to a wedge in the direction
of each nontrivial value of the Levi form of $M$ 
as was recently observed by {\sc Eastwood-Graham} \cite{EG99,EG01}
already in the case when $M$ is a hypersurface 
(see \cite[Example~1.2]{EG99} and Example~\ref{main-ex} below).
In this paper we propose an invariant geometric way of selecting 
those Levi form directions that are responsible for the extension, 
in particular, explaining the phenomenon discovered in \cite{EG99,EG01}.

The (extrinsic) Levi form of $M$ at a point $p$ can be seen
as  a hermitian form $L=L_p\colon T^c_pM\times T^c_pM\to (T_p\C^N/T_pM)\otimes\C$,
where $T^c_pM:=T_pM\cap iT_pM$ is the complex tangent space.
For every tangent vector $w\in T^c_pM$ consider the 
opening angle $0\le\g_w \le 2\pi$ of the plane convex cone $\C w \cap C_pV$,
where $C_pV$ denotes the closed tangent cone to $V$ at $p\in\d V$ (see \S\ref{basic}).
(We put $\g_w=0$ if $\C w \cap C_pV$ contains no open cone
and $\g_w=2\pi$ if $\C w \cap C_pV=\C w$.)
Define now the {\em Levi $(\pi/2)$-cone} of $V$ at $p$ as the convex hull
\begin{equation}\Label{maincone}
\G_p^{\pi/2} = \G_p^{\pi/2}(V) := {\sf conv}\,\{ L(w,w) : \g_w > \pi/2\}
\end{equation}
taken in the space $T_p\C^N/T_pM$.
We refer to $\g_w>\pi/2$ as the {\em angle condition}.
Throughout this paper we write 
$${\sf pr}\colon T_p\C^N \to T_p\C^N /T_pM$$ 
for the canonical projection.
Our main result provides holomorphic extension of CR-functions on $V$
in all additional directions of $\G_p^{\pi/2}(V)$:

\bt\Label{main}
Let $V$ be a wedge with edge $E$ at $p$ in a $C^4$-smooth submanifold  $M\subset\C^N$,
where $E$ is a $C^2$-smooth generic submanifold of $\C^N$.
Suppose that the interior of the cone $\Gamma_p^{\pi/2}\subset T_p\C^N/T_pM$ is nonempty. 
Then all continuous CR-functions on $V$ extend holomorphically to 
a wedge $V'$ in $\C^N$ with edge $E$ at $p$ such that 
\begin{enumerate}
\item[(i)] $\1{V'}$ contains a neighborhood of $p$ in $V$;
\item[(ii)] ${\sf pr}(C_pV') \supset \G_p^{\pi/2}$. 
\end{enumerate}
Moreover, the extension of each CR-function continuous up to $E$
is also continuous up to $E\cup V$.
\et

In particular, the wedge $V'$ is contained in the polynomial, rational and holomorphic 
(with respect to any Stein neighborhood) hulls of $V$.
The following example shows that $\pi/2$ for the cone $\G_p^{\pi/2}$ cannot be replaced by any smaller angle:
\be\Label{main-ex}
Consider the hypersurface
\begin{equation*}\Label{}
M:=\{z=(z_1,z_2,z_3)\in\C^3 : \Im z_3 = (\Im z_1)^2 - (\Im z_2)^2\}
\end{equation*}
and the wedge $V:=\{z\in M : \Im z_1 > |\Im z_2|\}$ with edge $E:=\R^3$ at $0$.
Then for $w := (1,1-\eps -i\sqrt{2\eps},0)\in T^c_0 M$ with $\eps>0$ arbitrarily small,
the Levi form $L(w,w)$ points down and the angle $\g_w<\pi/2$ is arbitrarily close to $\pi/2$.
However, any domain of holomoprhic extension of CR-functions on $V$
must be contained in the half-space $\Im z_3>0$.
Hence the Levi form of a vector $w$ with $\g_w$ arbitrarily close to $\pi/2$,
may not be in general a direction of extension as in Theorem~\ref{main}.
\ee

The assumption that the edge $E$ is {\em generic} in Theorem~\ref{main}
is also important and cannot be dropped as the following example shows:

\be
Consider the hypersurface
$M:=\{z=(z_1,z_2)\in\C^2 : \Im z_2 = |z_1|^2\}$
and, for a real number $1/2<\a<1$,
the wedge $V:=\{z\in M : \Im z_2 <\Im z_1^{1/\a} \}$ with edge $E:=\{z\in M : z_1=0\}$ at $0$.
Then, for $w:=(1,0)\in T^c_0M$, we have $L_0(w,w)>0$ and $\g_w=\pi\a>\pi/2$.
However, the CR-function $f(z_1,z_2):=(z_2-z_1^{1/\a})^{-1}$ on $V$
does not extend holomorphically to any wedge $V'$ as in Theorem~\ref{main}.
(It can be shown that the polynomial, rational and holomorphic hulls of $V$
do not contain any wedge in $\C^2$ with edge $E$, see e.g.\ \cite{ZZ2}.)
\ee

If one restricts attention to vectors $w\in T_p^cM\cap T_pE$,
then the angle condition in \eqref{maincone} is automatically satisfied.
Indeed, since $E$ is generic at $p$, it is easy to see that
the angle $\g_w$ is always either $2\pi$ or $\pi$ depending on whether $w$ belongs to $T_p^cE$ or not.
For Levi directions arising from vectors $w$ of this kind, a similar extension result is due
to {\sc Boggess-Polking} for $\g_w=2\pi$
and to {\sc Baracco} and the second author \cite{BZ01} for $\g_w=\pi$.
Another case has been treated by {\sc Eastwood-Graham} \cite{EG99,EG01},
where the cone $C_pV$ contains zero directions of the Levi form
with nonzero directions arbitrarily close. 
This situation can be also explained by Theorem~\ref{main},
since the angle $\g_w$ may be chosen arbitrarily close to $\pi$,
in particular, greater than $\pi/2$.
On the other hand, the following example shows that 
one may really need to consider $\g_w$ arbitrarily close to $\pi/2$
to get new directions of extension in $T_p\C^N/T_pM$:

\be\Label{noe}
Consider the hypersurface
\begin{equation*}\Label{}
M:=\{z=(z_1,z_2,z_3)\in\C^3 : \Im z_3 = |z_1|^2 + |z_2|^2 - (2+\eps)\Im z_1 \1z_2\}
\end{equation*}
for a small $\eps>0$. Then the Levi form is negative only for vectors $w\in T^c_0M$
whose direction is sufficiently close to that of $w_0=(-1+i,1+i,0)$,
for which $L_0(w_0,w_0)=-2\eps<0$. We choose $V:=\{z\in M : \Im z_1>0, \Im z_2>0\}$
so that $\g_{w_0}=\pi/2$.
Then there exist $w$ arbitrarily close to $w_0$ with $\g_w>\pi/2$ (see also \S\ref{equiv})
and hence Theorem~\ref{main} yields extension to a wedge on the negative
side of $M$.  On the other hand, for any vector $w$ with $L_0(w,w)<0$, the
angle $\g_w$ must be close to $\pi/2$.  \ee The technique of \cite{BZ01}
was based on attaching smooth ($C^{1,\b}$) analytic discs to $E\cup V$
which is possible in any direction $w\in T_p^cM\cap T_pE$ and then on
filling a wedge by these discs.  In the other directions $w\in T_p^cM$ with
$\pi/2<\g_w <\pi$ considered in this paper, such discs may not exist and
hence other methods are required.  Our approach here consists of two main
steps.  Using tools developed in \cite{ZZ2}, we first attach nonsmooth
analytic discs to $V$ that will not fill a wedge anymore but rather a
smaller region $V'$ that we call $\a$-wedge, where some directions are
bounded by $\a$th powers of others for some $\a>1/2$.  (Here $\a$ is chosen
such that $\g_w=\a\pi$.)  Then we attach analytic discs to properly chosen
submanifolds $\2M$ approximating $M$ such that, over a certain region,
$\2M$ is contained in $V'$.  For this step, the property $\a>1/2$ is
crucial and guarantees that the directions
 of $(1/\a)$th powers are not affected by 
the Taylor expansion of (the defining functions of) $M$. 

We conclude by giving an application of Theorem~\ref{main} to regularity
and holomorphic extension of CR-functions generalizing 
the ``edge-of-the-wedge'' theorem of {\sc Ajrapetyan-Henkin} \cite{AH81}.

\bc\Label{sum}
Let $M\subset\C^N$ be a $C^4$-smooth submanifold,
$V_1,\ldots,V_n$ be a wedges in $M$ with the same 
edge $E$ at $p$, where $E$ is a generic $C^2$-smooth submanifold.
Suppose that 
\begin{enumerate}
\item[(i)] $C_p V_1 + \cdots + C_p V_n = T_pM$;
\item[(ii)] $\G_p^{\pi/2}(V_1) + \cdots + \G_p^{\pi/2}(V_n) = T_p\C^N/T_pM$.
\end{enumerate}
Then every continuous function $f$ on $E$ that admits continuous CR-extension to each $V_j$,
also extends holomorphically to a neighborhood of $p$ in $\C^N$.
In particular, $f$ has the same smoothness as $E$
and its extension to each $V_j$ the same smoothness as $M$.
\ec

In the case of two wedges $V_1,V_2$ with opposite directions in a hypersurface $M$
containing zero directions of the Levi form (of $M$) having mixed signature
similar results are due to {\sc Eastwood-Graham} \cite{EG99,EG01}.

\section{Equivalent version of the angle condition}\Label{equiv}
The angle condition $\g_w>\pi/2$ has a simple interpretation
in terms of the directional cone of $V$.
For any generic submanifold $E\subset M\subset\C^N$ at a point $p\in E$ we have the canonical isomorphism
$$
(iT_pE\cap T_pM) / T^c_pE \to T_pM/T_pE.
$$
Then, if $V\subset M$ is a wedge in $M$ with edge $E$ at $p$,
the (closed) cone $C_pV$ can be represented, modulo $T^c_pE$, by an
open cone $i\Sigma\subset iT_pE\cap T_pM$ transversal to $T_pE$ such that 
the interior of $C_pV$ coincides with $T_pE+i\Sigma$.
We have:

\bl\Label{eq-cond}
For $w\in T^c_pM$, the condition $\g_w>\pi/2$ holds if and only if,
for suitable $e^{i\theta}\in S^1$, $\2w\in T^c_pM$ arbitrarily close to $w$
and for the decomposition 
$$e^{i\theta}\2w=a+ ib\in T_pE+iT_pE \mod T^c_pE,$$ 
one has $b\pm a\in\Sigma$.
\el

\bpf
By adding a suitable vector $v\in\Sigma$ to $w$,
we can assume that the plane $\C w$ intersects the interior of $C_pV$
in a plane cone of angle greater than $\pi/2$.
The last condition holds if and only if, for some choice of $e^{i\theta},a,b$ as above, one has
\begin{equation}\Label{ab}
\Im [(1-\tau)^\a (a+i b)]\in \Sigma \mod T^c_pE \text{ for all } \tau\in\Delta
\end{equation}
for some $\a>1/2$, where $\Delta\subset\C$ denotes the unit disc.
The inclusion in \eqref{ab} is clearly equivalent to
$b+ta\in \Sigma$ for all $|t|<1+\eps$ for some $\eps>0$
and hence to $b\pm a\in\Sigma$.
\epf

A particular case when $a=0$ that is $e^{i\theta}w=ib\in i\Sigma + T^c_pM$
or $w\in T^c_pM$ was considered in \cite{BZ01}.
It is clear that this happens precisely when either $\g_w=\pi$ or $\g_w=2\pi$.

\section{Preliminaries}\Label{basic}

We begin by defining open subsets in real manifolds having cone property.
By a {\em cone} $\G$ in the euclidean space $\R^m$ we always mean a subset 
invariant under multiplication by positive real numbers.
A subcone $\G'\subset\G$ is said to be {\em strictly smaller}
if $\1{\G'}\setminus\{0\}$ is contained in the interior of $\G$.

\bd\Label{cone-prop}
An open subset $V$ in a real $C^1$-smooth manifold $M$
is said to have the {\em cone property} at a boundary point $p\in\d V$
with respect to an open cone $\G\subset T_pM$ if, 
for any strictly smaller subcone $\G'\subset \G$
and some (and hence any) local coordinates on $M$ in a neighborhood of $p$,
one has $x+y\in V$ for all $x\in V$ and $y\in \G'$
sufficiently close to $p$ and $0$ respectively,
where the sum is taken with respect to the given coordinates.
\ed

We also say that $V$ has the cone property without specifying $\G$ if
there exists $\G$ for which the cone property holds.
It is clear from the definition that, if a given $V\subset M$
has the cone property with respect to each of two cones $\G_1$ and $\G_2$,
then the cone property also holds for the sum $\G_1+\G_2$.
Furthermore, among all such cones there is unique maximal one,
namely the sum of all of them that is automatically convex.
We call it the {\em tangent cone} to $V$ at $p$ and denote by $C_pV$.

An important special case of a subset with cone property is a wedge. 
We call an open subset $V\subset M$ a {\em wedge with edge $E$} at $p$,
where $E$ is a submanifold of $M$ through $p$,
if $V$ has cone property at $p\in\d V$ and if $E\subset\d V$.
In this case the tangent cone $C_pV$ is automatically invariant
under addition of vectors in $T_pE$ and 
can be canonically represented by a normal cone $\Sigma\subset T_pM/T_pE$
that we call the {\em directional cone} of the wedge $V$ at $p$. 

Finally, recall that a real submanifold $E\subset\C^N$
is called {\em generic at a point $p\in E$} if $T_pE + iT_pE = \C^N$.
Generic submanifolds are precisely the noncharacteristic submanifolds
for the system of Cauchy-Riemann equations.
If $M\subset\C^N$ is another submanifold containing $E$
which is generic at $p$, then $M$ is also generic at $p$. 
As usual continuous functions are called CR-functions
if they are solutions of the tangential Cauchy-Riemann equations
in distributional sense (see e.g.\ \cite{BERbook}).

\section{Attaching nonsmooth analytic discs}\Label{s2}

We discuss now about attaching analytic discs to a generic
submanifold $M$  of $\C^N$. The reader is referred 
to \cite{B91,T96,BERbook} for details. 
By an {\em analytic disc} in $\C^N$ we mean a holomorphic mapping  $A$ from 
$\Delta$, the unit disc in $\C$, into $\C^N$ which extends
at least continuously to $\1\Delta$; 
we still denote by $ A$ the image $A(\Delta)$.
We say that $A$ is {\em attached} to $M$ if 
$A(\partial \Delta) \subset M$. Usually one considers smooth 
($C^{k,\beta}$ with $k\geq 1$) discs. (The reason of the fractional 
regularity is due to the continuity of the Hilbert transform in the 
corresponding classes of functions.) We will also make an extensive 
use of the main technical tools developed in \cite{ZZ2}  
for attaching analytic discs that are $C^1$ except at $\tau=1$, 
where they are only H\"older-continuous.  For 
$\frac12<\alpha<1$, we put 
$$\beta:=2\alpha-1.$$ 
For any real $\delta$, 
consider the principal branch of $(1-\tau)^\delta$ on $\Delta$
that is real on the interval $(0,1)$; 
note that $(1-\tau)^{2\alpha}\in C^{1,\beta}$. 
As in \cite{ZZ2}, define $\Cal P^\alpha(\partial\Delta)\subset C^{\alpha}(\d \Delta)$ to be 
the subspace of all linear combinations of $(1-\tau)^\alpha$ and functions in $C^{1,\beta}(\partial\Delta)$
and $\Cal P^\alpha(\bar\Delta)\subset C^{\alpha}(\1\Delta)$ 
to be the subspace of linear combinations of $(1-\tau)^\alpha$, $(1-\bar\tau)^\alpha$ 
and functions in $C^{1,\beta}(\1\Delta)$.
These subspaces become Banach spaces
with the norm of a disc being the sum of the $C^\a$-norm of its nonsmooth factor
and the $C^{1,\beta}$-norm of the remainder. 
It is easy to see that the Hilbert transform 
is a bounded operator on $\Cal P^\alpha(\1\Delta)$. 
We choose coordinates $(z,w)\in\C^l\times\C^n=\C^N$, $z=x+iy$,
where $M$ is defined by $y=h(x,w)$ for a vector valued function $h=(h_1,\ldots,h_l)$ 
with $h(0)=0$, $h'(0)=0$. It is well-known (cf. \cite{B91})
that, given a sufficiently small holomorphic function $w(\tau)$ in $\Delta$ that is $C^\alpha$ 
in $\1\Delta$ and a sufficiently small vector $x\in\R^l$, there exists a (unique) small analytic disc 
$A(\tau)=(z(\tau),w(\tau))$ attached to $M$ in the same class with prescribed small value $x(1)=x$. 
It is shown in \cite{ZZ2} that the same statement also holds in the class $\Cal P^\alpha$
defined above and that the obtained disc $A$ dependens smoothly on $w(\cdot)$ and $x$. 
More precisely, one has:

\bp[\cite{ZZ2}]\Label{p8} Let $M$ be $C^{k+2}$-smooth ($k\geq 1$) and 
fix $\alpha> 1/2$. Then for every sufficiently small $x\in\R^l$ and 
holomorphic function $w(\cdot)\in\Cal P^\alpha(\1\Delta)$ there is an unique 
sufficiently small analytic disc $A=(z(\cdot),w(\cdot))\in\Cal P^\alpha$ with $x(1)=x$
attached to $M$ that depends in a $C^k$ fashion
on $w\in\Cal P^\alpha$ and $x\in\R^l$.
\ep

It is further shown in \cite{ZZ2} that, if $A(1)=0$
and if the coordinates $(z,w)$ are chosen as above,
the component $z(\tau)$ is in fact in $C^{1,\beta}$.

\section{Filling $\a$-wedges by nonsmooth analytic discs}

Our fist main step in proving Theorem~\ref{main}
will be to obtain a weaker statement
about the extension to a region smaller than a wedge.
More precisely, given submanifolds $M\subset M'\subset\C^N$
and an open subset, define an {\em $\a$-wedge in $M'$ over $V$
at a point $p\in \d V$ in the additional directions of an open convex cone $\G\subset T_pM'/T_pM$}
to be the open subset of $M'$ of the form
\begin{equation}\Label{d3}
V' = \{z\in W : \dist(z,V)< c\,\dist(z,\partial V)^{1/\a}\},
\end{equation}
where $W$ is a wedge in $M'$ with edge $M$ at $p$
and $c$ is a positive constant.
We say that an $\a$-wedge over $V$ is of class $C^{1,\delta}$
if the submanifold $M'$ can be chosen to be of class $C^{1,\delta}$.

The following statement is the main technical tool 
for the proof of Theorem~\ref{main}:

\bt\Label{t9}
Let $M\subset\C^N$ be a generic $C^4$-smooth submanifold  and 
$V\subset M$ be a domain with cone property at a point $p\in\partial V$.
Fix $\alpha> 1/2$ and $w\in T^c_{p}M$ such that
\begin{equation}\Label{7}
\g_w > \alpha\pi, \quad L_p(w,w)\ne 0.
\end{equation}
Then  analytic discs in $\C^N$ attached to $V$ fill an $\a$-wedge over $V$ at $p$ of class $C^{1,\delta}$
for some $0<\delta<1$ with one-dimensional additional direction arbitrarily close to that of $L_p(w,w)$.
\et

Here the new $\a$-wedge is understood in a submanifold $M'$
with boundary $M$ at $p$.
When $V$ is a wedge with generic edge $E$, then the arguments of 
{\sc Baouendi-Treves} \cite{BT81}
can be used to show that
continuous CR-functions on a neighborhood of $p$ in $V$ 
can be uniformly approximated on compacta by holomorphic polynomials.
Hence, replacing $V$ by such a neighborhood, we obtain from Theorem~\ref{t9}:

\bc\Label{extend}
Let $V$ be a wedge in a $C^4$-smooth submanifold $M\subset\C^N$
with generic edge $E$ at a point $p$.
Then all continuous CR-functions on $V$ admit CR-extension
to an $\a$-wedge over $V$ at $p$ of class $C^{1,\delta}$
for some $0<\delta<1$  with one-dimensional
additional direction arbitrarily close to that of $L_p(w,w)$.
\ec

Indeed, given a sequence of holomorphic polynomials converging uniformly
to a continuous CR-function on $V$ on compacta,
it also converges uniformly on the union of discs by the maximum principle.
Hence, if the discs fill an $\a$-wedge over $V$, the sequence of polynomials
converges there to a continuous CR-function.

\br
If a CR-function is of class $C^\kappa$ for some $1\le\kappa\le\infty$, the sequence of polynomials
can be chosen to converge in the $C^\kappa$-norm and therefore the 
CR-extension is also $C^\kappa$.
\er

\bpf[Proof of Theorem~{\rm\ref{t9}}]
Choose coordinates $(z,w)\in\C^l\times \C^{n}, 
z=x+iy$, with $p=0$ in which $M$ is defined by $y=h(x,w)$, where 
$h:=(h_1,\ldots,h_l)$ has no harmonic terms of order 2. 
Let $w=w_0\in T^c_p M$ satisfy \eqref{7}. 
It follows from Lemma~\ref{eq-cond} that $w_0$ can be chosen 
such that $(1-\tau)^{\alpha'} w_0 \in\Gamma$ for some $\alpha'>\alpha$ and all $\tau\in \Delta$,
where $\Gamma$ is a suitable strictly smaller subcone of the interior of $C_pV$.
For a small real parameter $\eta>0$ set 
$$
w(\tau)=w_\eta(\tau):= \eta(1-\tau)^\alpha w_0.
$$
We attach discs $A(\tau)=A_\eta(\tau)$ to $M$ with ``$w$-componets" 
$w(\tau)$ (the above choice of $w_0$ will guarantee that our discs with be automatically attached to $V$). 
We write 
$$z(\tau)=u(\tau)+ i v(\tau), \quad \tau=re^{i\theta}\in\1\Delta,$$ 
and solve the Bishop's equation 
$$u(\tau) = -T_1 h(u(\tau),w(\tau)), \quad \tau=e^{i\theta}\in\partial\Delta,$$
where $T_1$ is the Hilbert transform on $\d\Delta$ normalized by the condition 
$T_1(\cdot)|_{\tau=1}=0$.
By Proposition~\ref{p8}, for small $\eta$, this equation has solution $u(\cdot)$ in $\Cal P^\alpha$ with $u(1)=0$.
It is proved in \cite{ZZ2} that $u(\cdot)$ is in $C^{1,\beta}$ and moreover 
$\eta\mapsto v_\eta$, $\R\to C^{1,\beta}$ is of class $C^k$. 
Then also $v:=T_1 u\in C^{1,\beta}$ and $\eta\mapsto v_\eta$ is of class $C^k$. 
In particular, for the radial derivative $\partial_rv$, we have
\begin{equation}\Label{not5}
\partial_rv = \partial_rv|_{\eta=0} + (\partial_r\dot v|_{\eta=0})\eta + (\partial_r \ddot v|_{\eta=0})\eta^2 + o(\eta^2),
\quad \eta\to 0,
\end{equation}
where the dots stand for the derivatives in $\eta$.
Since $h(0)=0$ and $h'(0)=0$, we have $ v|_{\eta=0}\equiv0$ and $\dot v|_{\eta=0}\equiv0$.
Hence also $\dot u|_{\eta=0}\equiv0$ because 
$\dot u$ is related to $\dot v$ by the Hilbert transform. As for 
$\ddot v$, note that double differentiation in $\eta$ 
of $v(\cdot)=h(u(\cdot),w(\cdot))$ along $\partial\Delta$ yields
\begin{equation}\Label{not6}
\ddot v|_{\eta=0}= L_p(w_0,w_0)|1-\tau|^{2\alpha},\quad \tau\in\partial\Delta.
\end{equation}
Note that for $\tau=e^{i\theta}$, we have  
$|1-\tau|^{2\alpha}=|\theta|^{2\alpha}+o(|\theta|^{2\alpha})\in 
C^{1,\beta}(\d\Delta)$. Hence the harmonic extension of $|1-\tau|^{2\alpha}$ from $\partial\Delta$ 
to $\Delta$ is also of class $C^{1,\beta}(\d\Delta)$. 
By the Hopf Lemma, the radial derivative of this harmonic extension at $\tau=1$ is negative
and we have by \eqref{not6},
$$
\partial_r\ddot v|_{\tau=1,\eta=0} = c L_p(w_0,w_0), \quad c<0. 
$$
Hence  for sufficiently small $\eta$, we conclude from \eqref{not5} that $A_\eta=(z_\eta,w_\eta)$ 
is transversal to $M$ and $\d_r v$ points to a direction arbitrarily
close to $L_p(w_0,w_0)$. 
The conclusion about filling an $\a$-wedge over $V$ follows now from \cite{ZZ2}.
\epf

If \eqref{7} holds for several vectors $w_1,\ldots,w_l\in T^c_p M$
such that $L_p(w_1,w_1),\ldots,L_p(w_l,w_l)$ is a basis in $T_p\C^N/T_pM$,
then CR-functions on $V$ extend to a full-dimensional $\a$-wedge over $V$ 
whose additional directions can be chosen to be in any cone
strictly smaller than the convex span of the above basis. 
This result will follow from the proof of the ``edge-of-the-wedge'' theorem
of {\sc Ajrapetian} and {\sc Henkin} \cite{AH81} given by {\sc Tumanov} in \cite{T96}
using attached analytic discs (the original proof in \cite{AH81}
uses integral representations).
More precisely, for any positive real number $\a$, 
define the Levi $\a\pi$-cone $\G_p^{\a\pi}\subset T_p\C^N/T_pM$ 
to be the convex hull
$$
\G_p^{\alpha\pi} := {\sf conv}
\{L_p(w,w) :  w\in T^c_p M, \, \g_w > \alpha\pi\}.
$$

\bc\Label{c11}
In the setting of Theorem~{\rm\ref{t9}},
assume that $\G_p^{\alpha\pi}$ has a nonempty interior in $T_p\C^N/T_pM$.
Then for any strictly smaller subcone $\G'\subset \G_p^{\alpha\pi}$,
analytic discs attached to $V$ fill an $\a$-wedge $V'$ over $V$ at $p$
with additional directions of $\G'$.
In particular, all continuous CR-functions on $V$ 
that are approximated by polynomials on compacta
admit holomorphic extension to $V'$.
\ec

\br
If the boundary $\d V$ does not contain a generic submanifold $E$,
the classical arguments of \cite{BT81} cannot be applied to show
that any continuous CR-function on $V$ is approximated by polynomials
on compacta, even after replacing $V$ by a smaller neighborhood of $p$ in $V$.
Nevertheless it is shown in \cite{ZZ2} by a different construction
that holomorphic extension of CR-functions to an $\a$-wedge can be obtained also in this case.
\er

\bpf[Proof of Corollary~{\rm\ref{c11}}]
We write $\G_1\subset\subset \G_2$ if $\G_1$ is strictly smaller than $\G_2$.
Given $\G'\subset\subset \G_p^{\alpha\pi}$ as above, consider a 
polyhedral approximation of $\G'$ in $\G_p^{\alpha\pi}$:
\begin{equation}\Label{include}
\G'\subset \subset {\sf conv}\{ L_p(w_1,w_1),\dots, L_p(w_s,w_s)\}\subset \G_p^{\alpha\pi}.
\end{equation}
By Corollary~\ref{extend}, all (continuous) CR-functions on $V$ admit CR-extensions to 
$\a$-wedges $\2 V_1,\ldots,\2V_s$  over $V$ at $p$ of class $C^{1,\delta}$ with additional directions 
of vectors $v_1,\ldots,v_s$ arbitrarily close to those of $L_p(w_1,w_1),\dots, L_p(w_s,w_s)$ respectively.
Denote by $W_1,\ldots,W_s$ the corresponding wedges with edge $M$ at $p$ as in \eqref{d3}.
By \eqref{include}, we can assume 
$\G'\subset\subset {\sf conv}\{v_1,\dots,v_s\}$.
Then there exists a constant $C>0$ such that for any $z_0\in V$ and $1\le j\le s$,
\begin{equation}\Label{incl}
\{z\in W_j : \|z-z_0\| < c\, \dist(z_0,\d V)^{1/\a}\} \subset \2V_j.
\end{equation}
By repeating the arguments of the proof of \cite[Theorem~4.1]{T96}
one can show that all continuous CR-functions on $V$ 
admit holomorphic extensions to a subset
\begin{equation}\Label{subset}
\{z\in B(p) + \G' : \|z-z_0\| < c'\, \dist(z_0,\d V)^{1/\a}\},
\end{equation}
where $B(p)$ is a neighborhood of $p$ in $M$ and $c'>0$ is a constant,
both chosen independently of $z_0$
(we have identified $\G'$ with a cone in $\C^N$ in the normal direction to $T_pM$).
Such a choice can be obtained by rescaling the submanifold in \eqref{incl}
and using the uniformity of the construction of {\sc Tumanov} \cite{T96}
with respect to the $C^{1,\delta}$-norms of the submanifolds
and using his regularity results \cite{T93} for parameter dependence of solutions of Bishop's equation.
It follows now from \eqref{d3} that the subset \eqref{subset}
contains an $\alpha$-wedge over $V$ with additional directions of $\G'$.
\epf

\section{Passing from $\a$-wedges to wedges}

In this section we will assume that $V\subset M$ is a wedge with generic edge $E$
at $p$. We denote by $(C_pV)^0$ the interior in $T_pM$ of the tangent cone $C_pV$ of $V$ at $p$.
As before we fix $\a>1/2$.

\bp\Label{t12}
Let $V$ be a wedge in $M$ with generic edge $E$ at $p$
and $V'$ be an $\a$-wedge over $V$ at $p$
with one-dimensional additional direction.
Then the union of analytic discs attached to $V\cup V'$ contains
a wedge $V''$ with edge $E$ at $p$ in a submanifold $M''\subset\C^N$ with $\dim M'' = \dim M+1$
such that $\1{V''}$ contains a neighborhood of $p$ in $V$
and ${\sf pr}(C_pV'')$ is arbitrarily close to the additional direction of $V'$.
In particular, all continuous functions on $V$ 
admitting CR-extension to $V'$ are also CR-extendible to a wedge 
with edge $E$ at $p$ as above.
\ep

\bpf
Since $E$ is generic,
by shrinking $V$ and $V'$ if necessary we may assume
that continuous functions on $V\cup V'$ that are CR on $V'$
are uniformly approximated by holomorphic polynomials on compacta.
By choosing suitable holomorphic coordinates $z=x+iy\in\C^N$ in a neighborhood of $p$
we may assume that $p=0$,
$$T_0 E\subset \{y_{N-1}=y_N=0 \}, \quad  T_0 M\subset \{y_N=0\}, \quad (0,\ldots,0,i,0)\in (C_pV)^0$$
and the additional direction of $V'$ is that of $(0,\ldots,0,i)$.
We write $z=(z'',z_{N-1},z_N)\in\C^N$, $z'=(z'',z_{N-1})\in\C^{N-1}$
and denote by $\sigma\colon B(p)\to M$ any transversal projection from a neighborhood of $p$ in $\C^N$ to $M$.
In the sequel we write $c_j$ ($j=1,2,\ldots$) for positive constants
that may change during the proof.
Since $M$ is of class $C^2$, we have 
$$y_N>-c_1\|(z',x_N)\|^2, \quad z=(z',x_N+iy_N)\in M$$
with $c_1$ suitably chosen.
Let $W$ be a wedge with edge $M$ at $p$ as in \eqref{d3}.
% so that the additional direction of $V'$
%is given by ${\sf pr}(C_p W)\subset T_p\C^N/T_pM$.
Then there exist $c_1,c_2>0$ such that
\begin{equation}\Label{est1}
\{ z \in W : \sigma(z)\in V, \; c_1  \|(z',x_N)\|^2< y_N < -c_1\|(z',x_N)\|^2 + c_2\dist(z,\partial V)^{1/\a}\} \subset V'.
\end{equation}
Let further $\2V\subset V$ be a wedge with edge $E$ at $p$
whose directional cone at $p$ is a sufficiently small cone around $[(0,\ldots,0,1,0)]\in T_p M/T_p E$.
Then for $z$ near $p$ with $\sigma(z)\in \2V$, we have
\begin{equation}\Label{est2}
\dist(z,\partial V) > c_1|y_{N-1}|-c_2\|(z'',x_{N-1},x_N)\|^2
\end{equation}
for suitable constants $c_1,c_2>0$,
where we have used that $E$ is of class $C^2$. 
From \eqref{est1} and \eqref{est2} we conclude
\begin{equation}\Label{reg1}
\{z\in W : \sigma(z)\in \2V, \; c_1\|(z',x_N)\|^2 < y_N < - c_1\|(z'',x_{N-1},x_N)\|^2 + c_2 |y_{N-1}|^{1/\alpha}\}
\subset V'.
\end{equation}
We shall now define a deformation $\2M$ of $M$ 
whose part will enter the region in the left-hand side of \eqref{reg1}.
It is easy to see from the construction that, if $c_3$ is sufficiently large,
$$\2M:=\{ z\in \1W : y_N = c_3\|(z',x_N)\|^2 \}$$
is a submanifold of $\C^N$ near $p$ of the same dimension as $M$.
Futhermore we can shrink $\2M$ around $p$
and choose a sufficiently large constant $c_4$ for which \eqref{reg1} implies the inclusion
\begin{equation}\Label{inc4}
\{z\in \2M : \sigma(z)\in \2V,\;  y_{N-1}>c_4\|(z'',x_{N-1},x_N)\|^{2\alpha} \}\subset V'.
\end{equation}

We now construct an analytic disc attached to $\2M\cap V'$ that passes through $0$ and is transversal to $M$ there.
Consider a holomorphic function $z_{N-1}(\cdot)$ on the unit disc $\Delta$ 
that is $C^{1,\b}$ ($\b=2\a-1$) up to the boundary with $z_{N-1}(1)=0$ and satisfying
\begin{equation}\Label{t17}
y_{N-1}(\tau) \geq c_3|x_{N-1}(\tau)|^{2\alpha},\quad \tau\in\1\Delta,
\end{equation} 
where $c_3$ can be taken to satisfy $c_3 \ge 8c_4$.
By solving the Bishop's equation for sufficiently small $z_{N-1}$ (in the $C^{1,\b}$ norm), 
we can find a small analytic disc $z(\tau)$ with $z(1)=0$
and the nontrivial tangential derivative $\d_\theta z(1)\ne 0$, attached to $\2M$
with the component $z_{N-1}(\tau)$ prescribed above
such that the derivatives $\d_\theta z(\tau)$, $\tau\in\Delta$, are arbitrarily close to
the plane $\{(0,\ldots,0)\}\times\C\times \{0\}\subset T_pM$.
Then we can achieve $\sigma(z(\tau))\in \2V$ for all $z\in\Delta$.
Furthermore, since $\|(z''(\tau),z_N(\tau))\|\le c_5|1-\tau|$ with $c_5$ arbitrarily small
and $|x_{N-1}(\tau)| \ge c_6|1-\tau|$ for some $c_6$
and all $\tau\in\d\Delta$ in a suitable neighborhood $U$ of $1$, we conclude
$$
|x_{N-1}(\tau)|^{2\alpha}\geq \frac12|x_{N-1}(\tau)|^{2\alpha} + c_7\|(z''(\tau),x_N(\tau))\|^{2\alpha}
$$
with $c_7$ abitrarily large.
Therefore $z=z(\tau)$ is contained in the left-hand side of \eqref{t17} for $z\in \d\Delta\cap U\setminus\{1\}$.
On the other hand, it is clear from the construction that
$z(\tau)$ can be chosen to be contained in the left-hand side of \eqref{t17}
also for $\tau\in \d\Delta\setminus U$. It follows that $z(\tau)\in V'$
for all $\tau\in \d\Delta\setminus \{1\}$.
Furthermore, due to the construction of $\2M$ and the Hopf lemma, 
the radial derivative $\d_r z(1)$ 
is arbitrarily close to the subspace $\{(0,\ldots,0)\}\times\C^2\subset\C^{N-2}\times\C^2$
and $\d_r z(1)\notin T_p M$. 
Hence sufficiently small deformations $\2z(\cdot)$ of $z(\cdot)$
with $\2z(1)\in V$ fill a wedge $V''$
satisfying the required conclusion. 
\epf

\bpf[Proof of Theorem~{\rm\ref{main}}]
Fix $w\in T^c_pM$ with $\g_w>\pi/2$.
We may assume that continuous CR-functions on $V$
are uniformly approximated by holomorphic polynomials on compacta.
By Corollary~\ref{extend}, all continuous CR-functions on $V$
extend holomoprhically to an $\a$-wedge $V'$ over $V$ at $p$
with one-dimensional additional direction arbitrarily close to that of $L_p(w,w)$.
Then we conclude from Proposition~\ref{t12}
that CR-functions on $V$ also extend to wedge $V''$
in a submanifold $M''$ with $\dim M'' = \dim M+1$
satisfying the conclusion of the proposition.
In this way we obtain CR-extension to a $(\dim M+1)$-dimensional wedge $V''$
with ${\sf pr}(C_pV'')$ arbitrarily close to any given Levi form value $L_p(w,w)$
for $\g_w>\pi/2$.

If a CR-function on $V$ is continuous up to $E$,
it is approximated by polynomials uniformly up to $E$
and hence the extension is continuous up to $E\cup V$.
For these function, their holomorphic extension to a full-dimensional wedge $V'$
with ${\sf pr}(C_pV')\supset \G_p^{\pi/2}$ can be obtained
from the mentioned above proof of a strong version of the edge-of-the-wedge theorem 
given by {\sc Tumanov} \cite{T93,T96}.

If a CR-function is not known to be continuous up to $E$,
the edge $E$ can be ``pushed'' inside $V$ reducing
the extension problem to the above situation.
Then the conclusion follows from the parameter version of the edge-of-the-wedge theorem
following the arguments of \cite{T93,T96}.
\epf

\bpf[Proof of Corollary~{\rm\ref{sum}}]
If each of the cones $\G_p^{\pi/2}(V_1),\ldots,\G_p^{\pi/2}(V_n)$ has nonempty interior,
the conclusion follows directly from Theorem~\ref{main} and from 
the edge-of-the-wedge theorem  \cite{AH81}.
In general, the proof uses Corollary~\ref{extend} and Proposition~\ref{t12} as above
to obtain CR-extension to lower-dimensional wedges and 
then again the edge-of-the-wedge theorem. For the last passage
one needs the stronger version of the edge-of-the-wedge theorem assuming only $C^{1,\delta}$-regularity
of the submanifolds that follows from \cite{T93,T96}.
\epf

\end{document}